\documentclass[reqno]{article}
\usepackage{amsmath,amsthm,amssymb}
\usepackage[utf8x]{inputenc}
\usepackage[english]{babel}
\usepackage{amsmath}
\usepackage{amsfonts}
\usepackage{amssymb}
\usepackage{latexsym}
\usepackage{graphicx}
\usepackage{psfrag}
\usepackage{longtable}
\usepackage{mathrsfs}
\usepackage{enumitem}
\usepackage{cancel}
\usepackage{array}

\newtheorem{theorem}{Theorem}
\newtheorem{lemma}{Lemma}
\newtheorem{claim}{Claim}

\newtheorem{prop}{Proposition}
\newtheorem{rema}{Remark}

\newtheorem{property}{Property}
\newtheorem{question}{Question}

\begin{document}

\title{On more variants of the Majority Problem}

\author{ Paul-Elliot Angl\`es d'Auriac \thanks{%
University Paris-Est Cr\'eteil, LACL, 61 avenue du G\'en\'eral de Gaulle, F-94015 Cr\'eteil Cedex, France ; 
panglesd@lacl.fr
}
\and%
Francis Maisonneuve \thanks{%
Mines-ParisTech, 60 boulevard Saint-Michel 75272 Paris, France ; francis.maisonneuve@gmail.com}
\and%
Vivien Maisonneuve \thanks{%
v.maisonneuve@gmail.com}
\and%
Emmanuel Preissmann \thanks{%
emmanuel.preissmann@gmail.com}
\and%
Myriam Preissmann \thanks{
Univ.~Grenoble Alpes, CNRS, G-SCOP, Grenoble, France ; myriam.preissmann@grenoble-inp.fr}}

\maketitle

\begin{abstract}
The problem we are considering is the following. A colorblind player is given a set $\mathcal B=\{b_1,b_2, \ldots, b_N\}$ of $N$ colored balls. He knows that  each ball is colored either red or green, and that there are less  green than red balls (this will be called a {\it Red-green coloring}), but he cannot distinguish the two colors.
For any two balls he can ask whether they are colored the same. His goal is to determine the set of all green balls of $\mathcal B$ (and hence the set of all red balls). We study here the case where the Red-green coloring is such that there are at most $p$ green balls, where $p$ is given ; we denote by $Q(N,p, \le)$ the minimum integer $k$ such that there exists a method that finds for sure, for any Red-green coloring, the color of each ball of $\mathcal B$ after at most $k$ (color) comparisons. We extend the cases for which the exact value of $Q(N,p, \le)$ is known and provide lower and upper bounds for $Q(N,p, =)$ (defined similarly as $Q(N,p, \le)$, but for a Red-green coloring with exactly $p$ green balls).
\end{abstract}

\section{Introduction} 
There are several problems of the kind of the one considered in this paper. In the so-called "Majority problem", the goal is to find a ball of the "majority color" (red in our setting) within a minimum of comparisons. Saks and Werman \cite{SakWer91} and later Alonso et al \cite{AloReiSch93} and Wiener  \cite{Wie} have shown, with different methods, that $N - b(N)$ comparisons are necessary and sufficient to find for sure a red ball whatever the Red-green coloring, 
$b(N)$ being the number of $1$s in the binary representation of $N$.

To show the sufficiency it is enough to give a method using at most $N - b(N)$ comparisons. Since more complicated versions of this method and several parts of its validity proof will be used later,  we describe it in details (similarly as in  \cite {Aig}). The method uses boxes in which we place the balls, the {\it cardinality} of a box is the number of balls that it contains. During the procedure the contents of the boxes will vary : a box is called {\it monocolored} if all balls it contains are colored the same and else it is said {\it bicolored}. A bicolored box is said {\it balanced} if it has an equal number of balls of each color.  In the beginning we have $N$ boxes and we put exactly one ball in each, so at this stage all boxes are monocolored. At each step we compare two balls that belong to two distinct monocolored boxes of the same cardinality and then we merge the balls into one box. This box will stay monocolored if the two compared balls have the same color, else it will become balanced. The empty boxes are thrown away. 
So, at any time of this process the cardinality of each box is a power of $2$ and all bicolored boxes are balanced. We have to stop when no two monocolored boxes have the same cardinality. We claim that then all balls in the monocolored box of largest cardinality are red and we will now explain why.

As there is an equal number of red and green balls among the bicolored boxes, and since we have a Red-green coloring, there should be more red balls than green balls among the monocolored boxes. Let $2^k$ be the largest cardinality of a monocolored box. By the principle of the procedure, the number of balls that are in the other monocolored boxes is at most $\Sigma_{i=0}^{k-1}2^i$ (which value is set to $0$ in case $k=0$).  But $2^k >  \Sigma_{i=0}^{k-1} 2^i$, so the only way to have a majority of red balls among the monocolored boxes is that the balls in the box of size $2^k$ are red.

It remains to estimate the number of comparisons used by the method. Let $B_1, B_2, \ldots, B_r$ and $B_{r+1}, B_{r+2}, \ldots, B_{r+s}$ be respectively all monocolored and all balanced boxes obtained at the end of the procedure (where $r\ge 1$ and $s \ge 0$) and let $2^{k_i}$ be the cardinality of each $B_i$. We notice that the number of boxes decreases by $1$ after each comparison ; so the fact that the process ends with $r+s$ boxes means that we have done $N-(r+s)$ comparisons. Now, since each ball is in exactly one box, one has $N = \Sigma_{i=1}^{r+s} 2^{k_i}$ and this implies that $r+s \ge b(N)$. We can now conclude that we can provide one red ball after at most $N-b(N)$ comparisons.

\section{The Identification Problem}
We describe now the problems which are the subject of this paper. There is again a given set of $N$ balls colored with a Red-green coloring, but the goal now is to determine for sure the color of each of the balls within a minimum number of comparisons. Furthermore we will also consider more restricted kinds of Red-green colorings : given an integer $p$, a Red-green coloring of $N$ balls will be said {\it $p$-majored} if there are at most $p$ green balls and it will be said {\it $p$-equal} if there are exactly $p$ green balls.

Given two integers $N$ and $p<\frac{N}{2}$ we call {\it $(N,p, \le)$-identification} (respectively {\it $(N,p, =)$-identification}) the problem of determining for sure all the colors of $N$ balls colored by a $p$-majored Red-green coloring (respectively by a $p$-equal Red-green coloring) and we denote by $Q(N,p, \le)$ (respectively $Q(N,p, =)$) the minimum number of comparisons that are necessary to solve any instance of the $(N,p, \le)$-identification problem (respectively of the $(N,p, =)$-identification problem).

In order to describe algorithms solving the  $(N,p, \le)$-identification problem, or the $(N,p, =)$-identification problem, we use a method similar to the one above, using boxes. However, as we are now interested by the colors of all balls, we have to keep balls of distinct colors separated. For that purpose we will use {\it boxes divided into two sides} ; in the rest of the paper we will consider only boxes of this kind and will all along manage to have the balls in the same side of a box to be of the same color. 

Given a non empty box $B$ we will call {\it big side of $B$} the side of $B$ with the most balls. In case there is the same number of balls in each side, the big side of $B$ will be the one containing the ball $b_j$ where $j$ is the highest index of a ball in $B$. The side of $B$ which is not the big side will be called {\it the small side}. 

A box will be said of {\it Type $(x,y)$} if there are $x$ balls in its small side and $y$ in its big side (hence $x$ and $y$ are then two integers such that $x \le y$). 

A box whose small side is empty will be said {\it monocolored}, and a box of Type $(x,x)$ will be said {\it balanced} (the meaning of these definitions is the same as in the introduction). 

Initially we have no information on the colors and we put each ball $b_i$ in a monocolored box $B_i$ of Type $(0,1)$. All along we will manage to have a partition of the set $\mathcal B$ of the balls into boxes in such a way that the balls in the same side of a box are of the same color, different from the one of the balls in the other side. 
Then of course it has no sense to compare balls that are in the same box since the answer is already known.

Given two distinct boxes $B$ and $B'$ of a partition, the comparison of a ball $b$ in the big side of $B$ with a ball $b'$ in the big side of $B'$ leads to one of the two possible results : either the balls in the big sides of $B$ and $B'$ have the same color, or not. In any case, from the result of the comparison we may place all balls of $B$ and $B'$ into one box, in a way compatible with their colors.  Thus we obtain a new partition of the balls into boxes, expressing exactly all the knowledge on the balls colors provided by the comparison. If the Types of $B$ and $B'$ were $(x,y)$ and $(x',y')$ then the Type of the new box is either $(x+x',y+y')$ or $(min(x+y',x'+y),max(x+y',x'+y))$, according to whether the balls $b$ and $b'$ have the same color or not. 
For that reason we may consider that we do "boxes comparisons" rather than "balls comparisons". 

Notice that "doing a comparison" will always mean compare the colors of two balls {\bf and} merge the boxes containing these two balls as described above.

Then any sequence of $k$ comparisons results in $N-k$ non-empty boxes providing a partition $\mathcal P$ of the set $\mathcal B$ of balls into boxes. From $\mathcal P$  we can deduce the possible colorings of the balls :  a $p$-majored (resp. $p$-equal) Red-green coloring is {\it compatible} with $\mathcal P$  when any two balls in a same box are colored the same if and only if they are in the same side of the box. We will denote by {\boldmath $\mathcal C(\mathcal P)$} the set of $p$-majored (resp. $p$-equal) Red-green colorings that are compatible with a partition $\mathcal P$. Our goal is to reach, within a minimum number of boxes comparisons, a partition $\mathcal P$ such that $\vert \boldmath \mathcal C(\mathcal P)\vert = 1$.

Notice that the numbering of the boxes and the use of "small" and "big" sides are just a help for the description of the method. Also we notice that $\vert \boldmath \mathcal C(\mathcal P)\vert$ depends only on the number of balls in the sides of the boxes : a partition of the balls into boxes will be said {\it of Type} $(x_1, y_1)^{N_1} (x_2, y_2)^{N_2} \ldots (x_k, y_k)^{N_k}$ if it consists in exactly $N_i$ boxes of Type  $(x_i,y_i )$, for  $i=1$ to $k$.

\begin{rema} \label{rem} We notice that trivially {\bf $N-1$ comparisons are enough to solve both identification problems} since it will result in one box where only one side contains less than $p+1$ balls (and those are exactly the green balls of $\mathcal B$).
\end{rema}

\subsection{A tool and some useful Lemmas} 
 From our remarks above it is natural to represent a method $M$ solving the $(N,p, \le)$ or $(N,p, =)$-identification problem by a {\it labeled binary tree $T_M$}  as follows :
 \begin{itemize}
\item Each vertex $u$ is labeled with a partition $\mathcal P(u)$ into boxes,
\item The partition $\mathcal P(R)$ of the root $R$ of $T_M$ is made of $N$ boxes $B_1, B_2, \ldots, B_N$ where each $B_i$ contains only the ball $b_i$.
\item Each non-leaf vertex $u$ is furthermore labeled with  a couple $(B(u),B'(u))$ of boxes  belonging to $\mathcal P(u)$, and it has
\begin{enumerate} 

\item[-]  one child connected to $u$ by an edge labeled "=". This child is labeled by the partition obtained from $\mathcal P(u)$ by replacing $B(u)$ and $B'(u)$ by the new box resulting from the comparison of $B(u)$ or $B'(u)$ when the balls in the big sides of $B(u)$ and $B'(u)$ have the same color,

\item[-] one child  connected to $u$ by an edge labeled "$\neq$". This child is  labeled by the partition obtained from $\mathcal P(u)$ by replacing $B(u)$ and $B'(u)$ by he new box resulting from the comparison of $B(u)$ or $B'(u)$ when the balls in the big sides of $B(u)$ and $B'(u)$ have different colors.

\end{enumerate}

\item each  leaf-vertex $\ell$ is such that $\vert \boldmath \mathcal C(\mathcal P(\ell))\vert = 1$,

 \end{itemize}
 
 The binary tree $T_M$ is in a sense a "user manual" of the method $M$ as it gives at any step which boxes to compare. However, as it should contain all possible cases we will use more synthetic ways to describe a method. Nevertheless we will see in the following that $T_M$ may be useful to compute bounds on $Q(N,p,\le)$ and $Q(N,p,=)$.
 
 Notice that we can always assume that the two answers "$=$" or "$\neq$" may occur after a comparison since else this comparison is useless. So if $v$ and $w$ are the two children of a vertex $u$ of $T$ then $\mathcal C(\mathcal P(u))$ is equal to the disjoint union of the two non empty sets $\mathcal C(\mathcal P(v))$ and $\mathcal C(\mathcal P(w))$.

 A vertex of $T_M$ will be said at {\it level $k$} if it is at distance $k$ from the root. The {\it  height} of $T_M$ is the maximum level of a vertex of $T_M$.

Let $v$ be a vertex of $T_M$ at level $k \ge 1$. We call {\it parent of} $v$, denoted by $p(v),$ the unique neighbor of $v$  which is at level $k-1$. The vertex $v' \neq v$ such that $p(v')=p(v)$ will be said the 
 {\it brother of} $v.$ Remark that in a binary tree, any path $P$ issued from the root  contains $p(v)$ for each $v\in P$ which is not the root.

 From its definition and some preceding remarks we get the following properties of $T_M$.
 
 \begin{property} \label{propT}
 Any binary tree $T_M$ associated to a method $M$ solving the ($N,p, \le$)-identification problem satisfies the following properties:
\begin{itemize}
 \item[(i)] each $p$-majored Red-green coloring of $\mathcal B$ labels exactly one leaf,
 
 \item[(ii)] the partition of a vertex at level $k$ consists in $N-k$ boxes obtained after $k$ comparisons,
 
 \item[(iii)] the maximum number of comparisons used by an execution of $M$ is equal to the height of $T_M$.
 \end{itemize}
 \end{property}

Notice that similar remarks and properties are valid in the case of the $(N,p,=~)$-identification problem. 

From the (easy and wellknown) fact that a binary tree of height $k$ has at most $2^k$ leaves, Property \ref{propT}(i)  and Property \ref{propT}(iii), we get the following lower bound that we will use later.

\begin{lemma} \label{trouver}
For every two integers $N$ and $p<\frac{N}{2}$ we have :

$Q(N,p, =) \ge \log_2({{N} \choose {p}})$,

\end{lemma}

Let $\perp$ denotes either $"="$ or $"\le"$. The two following Lemmas are also useful.
\begin{lemma} \label{borneinfcomm} 
For every two integers $N$ and $1 \le p<\frac{N}{2}$ we have :
 $$Q(N,p,\perp) \ge 2 +Q(N-2, p-1, \perp).$$
 \end{lemma} 

\proof
Let  $T_M$ be the tree associated to a method solving the $(N,p,\perp)$-identifica-tion problem within $Q(N,p,\perp)$ comparisons. Since $p\ge 1$ the height of $T_M$, which is equal to $Q(N,p,\perp)$ by Property \ref{propT}(iii), is at least 2 and there is in $T_M$ an edge  labeled "$\neq$" between the root $R$ and a vertex $v$. In $\mathcal P(v)$, the first two balls $b,b'$ that are compared by $M$ are in a box $B^*$ of Type $(1,1)$ and all other boxes  are of Type $(0,1)$ and contain the balls of $\mathcal B \setminus\{b,b'\}$. Notice that the set $\mathcal B \setminus\{b,b'\}$ contains $N-2$ balls and exactly one less green ball than $\mathcal B$. On another hand, for any leaf $\ell$ of the subtree $T_v$ of $T_M$ of root $v$ we have $\vert  \mathcal C(\mathcal P(\ell))\vert=1$  and hence the balls $b$ and $b'$ should be in an unbalanced box of $\mathcal P(\ell)$. This means that, on any path of $T_M$ from $v$ to a leaf, there is a vertex where one of the two boxes that are compared is $B^*$ containing only $b$ and $b'$. 

We claim that $T_v$ gives a method to find the $(p-1)$-equal, or $(p-1)$-majored (depending on the value of $\perp$), Red-green coloring of $N-2$ balls : use the $N-2$ boxes of cardinality $1$ of $\mathcal P(v) \setminus B^*$, ignore $b$ and $b'$ in the partitions, do the box-comparisons and follow the edgess as indicated by the labelings, except when you reach a vertex of $T_v$ labeled with the "empty" box $B^*$  : then don't do any comparison and go to any child. Then, for any leaf $\ell$ of $T_v$  the resulting partition $\mathcal P_v(\ell)$ is the same as the one obtained from $\mathcal P(\ell)$ by withdrawing $b$ and $b'$ that are in different sides of one box. Hence $\vert \mathcal C(\mathcal P_v(\ell))\vert =1$, else we would contradict the fact that $\vert  \mathcal C(\mathcal P(\ell))\vert=1$. So the height of $T_v$ minus $1$ (we skipped exactly one comparison) is at least $Q(N-2, p-1, \perp).$ Since the height of $T_v$ is at most the height of $T_M$ minus $1$ we get that $Q(N,p,\perp)-2 \ge Q(N-2, p-1, \perp).$

\qed

 \begin{lemma} \label{bornesupcomm} 
 For every two integers $N$ and $p < \frac{N-1}{3}$  we have :
 $$Q(N,p,\perp) \le p +Q(N-(p+1), p, \perp).$$
\end{lemma} 

\proof In order to prove the upper-bound on $Q(N,p,\perp)$ it is sufficient to exhibit a method to solve the $(N,p,\perp)$-identification problem within at most  $p +Q(N-p-1, p, \perp)$ comparisons.
So let us assume that we have a set of $N$ balls $b_1, b_2, \ldots, b_N$ colored by a $p$-equal or $p$-majored (depending on the value of $\perp$) Red-green coloring. We may first compare $b_1$ to $b_2, b_3, ...b_i$ for some $2 \le i\le N$, thus obtaining a partition of Type $(x,y)^1 (0,1)^{N-i}$ for some integers $0 \le x \le y$ such that $x+y=i$. The rule will be to do so until $y=x+p+1$, which will happen since we must have one side containing at most $p$ balls and $N > 3p+1$. At this stage, we have done $2x+p$ comparisons, all balls in the side of cardinality $y=x+p+1> p$ of the only bicolored box should be red and those $x\le p$ in the other side are green :  it remains $N-(p+2x+1)$ other balls whose colors have to be determined and the only information we have about them is that there are at most $p-x$, or exactly $p-x$ in case $\perp$ means "=", green balls among them. Since $x\le p$ and we assumed $N>3p+1$ then $N-(p+2x+1) > 2(p-x)$ and we can determine the colors of the remaining balls using at most $Q(N-(p+2x+1), p-x, \perp)$ comparisons. So $Q(N,p,\perp) \le Max_{0\le x \le p} (p+2x + Q(N-(p+2x+1), p-x, \perp))$. This maximum is attained for $x=0$ since by Lemma \ref{borneinfcomm} we have for any $1 \le x \le p:$
\begin{align*}
{p+Q(N-(p+1),p,\perp)}  &\ge  p+2 +Q(N-(p+3), p-1, \perp) \\
   & \ge  \ldots \\
   & \ge p+2x+Q(N-(p+2x+1), p-x, \perp).
\end{align*}
\qed

\subsection{Bounds and exact values for $Q(N,p, \le)$}

The theorem below is due to Aigner \cite{Aig} but we did a different proof (that we find simpler).

\begin{theorem} \label{prop1} \cite{Aig}
Let $N$ and $p$ be  integers such that $0 \le p < \frac{N}{2}$, we have:
$$Q(N,p, \le) \le N+1- \left \lfloor \frac{N+1}{p+1} \right \rfloor.$$
\end{theorem}

\proof 
When $p=0$ : $Q(N,0, \le) =0= N+1- \lfloor \frac{N+1}{1} \rfloor,$ and the bound is correct. Let us consider the case when $p \ge 1$.

Assume first that $N<3p+2$. Then $\lfloor \frac{N+1}{p+1} \rfloor = 2$ and  $N+1- \lfloor \frac{N+1}{p+1} \rfloor= N-1$ which is always an upper bound of $Q(N,p,\le)$, as already noticed in Remark \ref{rem}.

\noindent Let us consider now the case where $N \ge 3p+2$.
We may use consecutively Lemma \ref{bornesupcomm} as long as we have a number of balls which is at least $3p+2$ :
$$Q(N,p, \le) \le p+Q(N-(p+1),p,\le) \le \ldots \le \ell p +Q(N-\ell(p+1), p, \le).$$ Then $\ell \ge 1$ is such that $$2p+1 \le N- \ell(p+1) < 3p+2 \Longleftrightarrow 2p+2 \le N+1- \ell(p+1) < 3p+3$$
$$\Longleftrightarrow N+1-2(p+1) \ge \ell(p+1) > N+1-3(p+1)$$
 $$\Longleftrightarrow \ell = \left \lfloor \frac{N+1}{p+1} \right \rfloor -2.$$

Then we have 
\begin{eqnarray*}
  Q(N,p, \le) \le \ell p +Q(N-\ell(p+1), p, \le) & \le & \ell p+ N-\ell(p+1)-1 \\
   & = & N - \ell -1= N +1- \left \lfloor \frac{N+1}{p+1} \right \rfloor .
\end{eqnarray*}

\qed

\bigskip

In his paper Aigner  showed that the upper bound in Theorem \ref{prop1} is near from the exact value of $Q(N,p, \le)$ (Theorem 4 in \cite{Aig}). Our next Theorem \ref{prop2} extends its results and partially answers a question of Wildon (Problem 8.1 in \cite{Wildon16}).

\begin{theorem} \label{prop2}
Let $N$ and $p$ be  integers such that $0 \le p < \frac{N}{2}$, we have:
\begin{itemize}
\item $N- \lfloor \frac{N}{p+1} \rfloor \le Q(N,p, \le) \le N+1-\lfloor \frac{N+1}{p+1} \rfloor $,
\item $Q(N,p, \le) = N +1- \lfloor \frac{N+1}{p+1} \rfloor$ when $N \equiv r$ $[p+1]$ and $r=p$ or $0 \le r \le \lfloor \frac{p}{2} \rfloor$.
\end{itemize}
\end{theorem}

\proof
The upper bound of $Q(N,p, \le)$ is the one of Theorem \ref{prop1} and it is equal to $Q(N,p,\le)$ when $p=0$. We prove the lower bound in case $p\ge 1$.

Let $M$ be any method that solves the ($N,p, \le$)-identification problem. The binary tree $T_M$ has a leaf $\ell$ which is connected to the root by a path whose all edges are labeled with "=". The partition $\mathcal P(\ell)$ contains then only monocolored boxes. 

\noindent In case one box $B$ of $\mathcal P(\ell)$ is of cardinality at most $p$,  there will be at least two $p$-majored Red-green colorings compatible with $\mathcal P(\ell)$ : one where all balls are  red, and one where all balls in $B$ are green and all other balls are red. This is in contradiction with the fact that $\ell$ is a leaf of $T_M$. 

\noindent So all boxes of $\mathcal P(\ell)$ are of cardinality at least $p+1$ ; hence $\mathcal P(\ell)$ contains at most $\lfloor \frac{N}{p+1} \rfloor $ boxes. These are obtain after at least $N-\lfloor \frac{N}{p+1} \rfloor $ comparisons. Since $M$ was chosen as any method solving the ($N,p, \le$)-identification problem, this implies that indeed $Q(N,p, \le) \ge N- \lfloor \frac{N}{p+1} \rfloor$. 
\medskip

Assume now that $N \equiv p$ $[p+1]$. Then  $N+1-\lfloor \frac{N+1}{p+1} \rfloor = N+1-(\lfloor \frac{N}{p+1} \rfloor+1)= N- \lfloor \frac{N}{p+1} \rfloor $, hence the two bounds are equal and the second statement of the theorem is proved.

\medskip

It remains to consider the case where $N \equiv r$ $[p+1]$ for $0 \le r \le \lfloor \frac{p}{2} \rfloor$.

\noindent So, let us consider two integers $N$ and $q$ such that $N \equiv r$ $[p+1]$ for some $0 \le r \le \lfloor \frac{p}{2} \rfloor$ and $q=\lfloor \frac{N+1}{p+1} \rfloor=\lfloor \frac{N}{p+1} \rfloor$. By our already proven bounds we have  $N- q \le Q(N,p, \le) \le N+1-q$. Assume that there exists a method $M^*$ that solves the $(N,p, \le)$-identification problem using only $N- q$ comparisons. By Property \ref{propT}(iii), the height of $T_{M^*}$ is equal to $N-q$ and the partition associated to a leaf of $T_{M^*}$ should then contain at least $q$ boxes. As in the first part of our proof we know that there exists a leaf $\ell^*$ of $T_{M^*}$ whose associated partition  $\mathcal P^*(\ell^*)$ is made of at most $q$, and hence  exactly $q$,  monocolored boxes. 

Let  $(R, u_1, u_2, \ldots, u_{N-q}=\ell^*)$ be the path from the root $R$ to $\ell^*$ in $T_{M^*}$ and let 
$$i= \min \{1\le j \le N-q-1\vert \mathcal P^*(u_{j+1}) \text{ contains } q \text{ boxes of cardinality at least } p+1\}$$ 
($i$ does exist since $N-q-1$ satisfies the requirement). 

Consider the vertex $v_{i+1}$ connected  to $u_i$ by an edge labeled $"\neq"$ (that is $v_{i+1}$ is the brother of  $u_{i+1}$). From our assumptions, the partition $\mathcal P^*(v_{i+1})$ consists into $q-1$  monocolored boxes  of cardinality at least $p+1$ , called {\it big monocolored boxes of $\mathcal P^*(v_{i+1})$}, one bicolored box $B$ of cardinality at least $p+1$ containing $x$ balls in one side and $y$ in the other where $1 \le x \le y \le p$, and $N-(i+1)-q$ monocolored boxes of cardinality at most $p$, called {\it small monocolored boxes of $\mathcal P^*(v_{i+1})$}.

Considering $\mathcal P^*(v_{i+1})$, we denote by $\mathcal A$ the set of balls belonging to 
 big monocolored boxes, 
$\mathcal B_x$ the set of the $x$ balls in the small side of $B$, $\mathcal B_y$ the set of the $y$ balls in the big side of $B$, and $\mathcal C$ the set of balls belonging to  small monocolored boxes. Notice that $\mathcal A$, $\mathcal B_x$, $\mathcal B_y$, $\mathcal C$ is a partition of the whole set $\mathcal B$ of balls.

From $v_{i+1}$ we will follow a path in $T_{M^*}$ using the following rule that will all along ensure partitions with exactly one bicolored box, containing $\mathcal B_x$ in one side and $\mathcal B_y$ in the other, and monocolored boxes containing balls from $\mathcal A\cup \mathcal C$.

{\bf Rule} : We start from the vertex $v_{i+1}$ which by definition satisfies the requirement. As long as we are on a non-leaf vertex $v$ 
we go down in the tree by repeating the following~:
\begin{itemize}
\item if $B(v)$ and $B'(v)$ (the two boxes of $\mathcal P^*(v)$ that will be compared at this step) are both monocolored  we follow the edge labeled $"="$ and get to the next vertex $v$,
\item else one of them, say $B(v)$, is the bicolored box of $\mathcal P^*(v)$. In that case, if in $\mathcal P^*(v)$ the balls in $B'(v)$ are all in $\mathcal C$ then we follow the edge labeled with the answer that merges the balls in $B'(v)$ with those of $\mathcal B_x$ and else we follow the edge with the answer that merges the balls in $B'(v)$ with those of $\mathcal B_y$ ; and get to the next vertex $v$.
\end{itemize}

The process stops when we reach a leaf-vertex $l$. By the rule, $\mathcal P^*(l)$ has monocolored boxes whose contents are included in $\mathcal A\cup \mathcal C$ and exactly one bicolored box $B_l$. Let $\mathcal B_{l,x}$ be the set of balls of $B_l$ in the side containing $\mathcal B_x$ and let $\mathcal B_{l,y}$ be the set of balls of $B_l$ in the side containing $\mathcal B_y$. By definition the balls of $\mathcal B_{l,x}$ that are not in $\mathcal B_x$ are in  $\mathcal C$ and those of $\mathcal B_{l,y}$ that are not in $\mathcal B_y$ are in $\mathcal A \cup \mathcal C$.

\begin{claim} \label{B1+C}
$\vert \mathcal B_x \vert + \vert \mathcal C \vert \le p$.
\end{claim}

\proof There are at least $(q-1)(p+1)$ balls in $\mathcal A$ and exactly $N=q(p+1)+r$ balls in $\mathcal B$, so we have at most $p+1+r$ balls in $\mathcal B_x \cup \mathcal B_y \cup \mathcal C$. Since  $\vert \mathcal B_x \cup \mathcal B_y \vert \ge p+1$, one has $\vert \mathcal B_x \cup \mathcal B_y \vert = p+1+s$, for some $0 \le s \le r \le \lfloor \frac{p}{2} \rfloor$. Then $0 \le \vert \mathcal C \vert \le N-(q(p+1)+s) = r-s$, and $\vert  \mathcal B_x \vert \le \lfloor \frac{p+1+s}{2} \rfloor \le \lceil  \frac{p+1+s}{2} \rceil \le \vert  \mathcal B_y \vert \le p$, so that 
$$\vert  \mathcal B_x \vert + \vert  \mathcal C \vert \le \left \lfloor \frac{p+1+s}{2} \right \rfloor +r - s = \left \lfloor \frac{p+1-s}{2} \right \rfloor +r \le \left \lfloor \frac{p+1}{2} \right \rfloor + \left \lfloor \frac{p}{2} \right \rfloor = p.$$
The claim is proved.

\begin{claim} \label{nomonoc}
No monocolored box of $\mathcal  P^*(l)$ is included in $ \mathcal C$.
\end{claim}

\proof Let us assume the contrary and let $ \mathcal C_1 \subseteq  \mathcal C$ be the set of balls of a monocolored box of $\mathcal P^*(l)$. By definition $\mathcal B_{l,x} \cup \mathcal C_1 \subseteq \mathcal B_x \cup \mathcal C$, so by Claim \ref{B1+C} we have $\vert \mathcal B_{l,x} \cup \mathcal C_1 \vert \le p$. Then there are at least two $p$-majored Red-green colorings compatible with $\mathcal P^*(l)$ : one colors in green only the vertices of $\mathcal B_{l,x}$ and the other colors in green the vertices of $\mathcal B_{l,x} \cup \mathcal C_1$. This contradicts the fact that $l$ is a leaf of $T_M^*$ and ends the proof of the claim.

\bigskip
A few remarks are now enough to end the proof of the theorem.
Let us remind that $\mathcal P^*(v_{i+1})$ contains $N-(i+1)$ boxes and among them $N-(i+1)-q$ are monocolored boxes included in $\mathcal C$. By Claim \ref{nomonoc} all these monocolored boxes disappear after doing all the comparisons indicated by the path from $v_{i+1}$ to $l$ obtained by applying the Rule. 
As we already know, each of these comparisons decreases the number of boxes by $1$ and it is straightforward to verify that, furthermore, the number of monocolored boxes included in $\mathcal C$ either stays the same or decreases by $1$ (only in case at least one of the two boxes that are compared is a monocolored box included in $\mathcal C$). So the number of boxes in $\mathcal P^*(l)$ is at most $q$. 
Morover, by our assumption on $M^*$, the height of $T_{M^*}$ and hence the level of $l$ is at most $N-q$, so that the number of boxes in $\mathcal P^*(l)$ is at least $q$.  From all these facts we deduce that each comparison done to get $\mathcal P^*(l)$ from $\mathcal P^*(v_{i+1})$ concerns at least one  monocolored box included in $\mathcal C$. So, as we followed the Rule, the bicolored box $B_l$ of $\mathcal P^*(l)$ is such that $\mathcal B_{l,x} \subseteq \mathcal B_x \cup \mathcal C$ and $\mathcal B_{l,y}= \mathcal B_y $.

Then, by Claim \ref{B1+C} and the fact that $\vert \mathcal B_{l,y}\vert = y \le p$, each side of $B_l$ contains at most $p$ balls and we have at least two $p$-majored Red-green colorings compatible with $\mathcal P(l)$ : everything colored red except the balls in one of the two sides of $B_l$. This contradicts the fact that $l$ is a leaf of $T_M$. 
\qed

\bigskip

From the preceding results we get the following properties of $Q(N,p,\le)$.

 \begin{property} \label{propQ}

\begin{itemize}

\item[(1)] Let $p$ be a nonnegative integer. 

The function : $N \mapsto Q(N,p,\le)$ ($N \ge 2p+1$) is nondecreasing and $1$-lipschitz, that is :

$$Q(N,p,\le) \le Q(N+1,p, \le) \le 1+ Q(N,p, \le)$$

\item[(2)] Let $N$ be a nonnegative integer.  

The function : $p \mapsto Q(N,p,\le)$ ($ p <  \frac{N}{2}$) is nondecreasing.

\item[(3)] $Q(N+1,p+1,\le) \ge 1+Q(N,p,\le)$ for any $N,p$ such that $N \ge 2p+2$.

\item[(4)] If $N$ and $p$ are such that $N \ge 2p+1$ and $Q(N,p, \le) < N+1 - \lfloor \frac{N+1}{p+1} \rfloor$ then for any $N'=N+k(p+1)$ ($k\ge 0$) we have $Q(N', p, \le) = N'-\lfloor \frac{N'}{p+1} \rfloor$.
\end{itemize}

 \end{property}
 
 \proof
 \begin{itemize}

\item[(1)] From the first statement of Theorem \ref{prop2} one has $$Q(N,p, \le) \le N+1- \left \lfloor \frac{N+1}{p+1} \right \rfloor \le Q(N+1, p, \le).$$

Let $M$ be a method determining a $p$-majored Red-green coloring of $N$ balls within at most $Q(N, p, \le)$ comparisons.  It is easy to derive  from $M$ a method determining a $p$-majored Red-green coloring of $N+1$ balls using at most one more comparison : put aside one of the $N+1$ balls, you have then a $p$-majored Red-green coloring of $N$ balls, determine their colors using $M$, if you get less that $p$ green balls one more comparison will be necessary to know the color of the ball that was put aside.

\item[(2)] Obvious since, when $p+1 < \frac{N}{2}$, any  $p$-majored Red-green coloring is a $(p+1)$-majored Red-green coloring.

\item[(3)]  Let $N \ge 2p+2$. From $(1)$ we know that $Q(N, p,\le) \le 1 + Q(N-1,p,\le)$ and from Lemma \ref{borneinfcomm} one has $Q(N+1,p+1,\le) \ge 2 +Q(N-1, p, \le).$ Hence $$Q(N+1,p+1,\le) \ge  2 +Q(N, p,\le) - 1= 1 +Q(N, p,\le).$$

\item[(4)] 
From Lemma \ref{bornesupcomm} we deduce $Q(N', p, \le) \le kp +Q(N,p, \le)$ for any $N, p, N'$ such that $N \ge 2p+1$ and $N'= N+k(p+1)$ for some $k \ge 0$. 

If $Q(N,p, \le) < N+1 - \lfloor \frac{N+1}{p+1} \rfloor$ then by Theorem \ref{prop2} we know that $Q(N,p, \le~) = N - \lfloor \frac{N}{p+1} \rfloor$ 
; so we have:$$Q(N', p, \le) \le kp +N - \left \lfloor \frac{N}{p+1} \right \rfloor= N'-k- \left \lfloor \frac{N}{p+1} \right \rfloor= N'- \left \lfloor \frac{N'}{p+1} \right \rfloor.$$

\end{itemize}
 \qed
 \bigskip
 
 Given two integers $N$ and $p$ such that $0 \le p < \frac{N}{2}$, let us define $Q_+(N,p,\le~):=N+1-\lfloor \frac{N+1}{p+1} \rfloor $ 
 and $Q_-(N,p,\le) := Q_+(N,p,\le)-1.$

Theorem \ref{prop2} shows that $Q(N,p,\le)$ is always equal to either $Q_+(N,p,\le)$ or $Q_-(N,p,\le)$ and that is is equal to $Q_+(N,p,\le)$ whenever $N \equiv r$ $[p+1]$ where $r=p$ or $0\le r \le \lfloor \frac{p}{2} \rfloor$. Thanks to Property \ref{propQ} we have furthermore the following interesting property.
 
 \begin{property} \label{propr}
 For every nonnegative integer $p$ exactly one of the two following statements is satisfied :
 \begin{enumerate}
\item For every $N \ge 2p+1$ we have $Q(N,p,\le)=Q_+(N,p,\le).$
\item There exist two integers $\frac{p+1}{2} \le r_p \le p-1$ and $N_{p} \ge 2p+1$ such that for every $N\ge N_{p}$ we have : $$Q(N,p, \le)=Q_-(N,p,\le) \mbox{ if and only if } N \equiv r \mbox{ } [p+1] \mbox{ for } r_p \le r \le p-1.$$ 
\end{enumerate} 
\end{property}
 
\proof Let us assume that the first statement is not satisfied. Then the set $\mathcal S=\{N \vert N\ge 2p+1 \mbox{ and } Q(N,p,\le)=Q_-(N,p,\le)\}$ is not empty.
Then we may choose $r_p$ as the minimum $r$ for which there exists $N \in \mathcal S$ such that  $N \equiv r$ $[p+1]$. Let now $N_p$ be the minimum $N \in \mathcal S$ such that $N \equiv r_p$ $[p+1].$ By Theorem \ref{prop2} and the definition of $r_p$ we know that $\frac{p+1}{2} \le r_p \le p-1$ and that $Q(N,p, \le) =Q_+(N,p,\le)$ for every $N \ge 2p+1$ such that $N \equiv r$ $[p+1]$ for $0\le r \le r_p-1$ or $r=p.$

Assume now that $Q(N,p,\le)=Q_-(N,p, \le)$ for some $N \equiv r$ $[p+1]$ where $r< p-1.$ By Theorem \ref{prop2} and our hypothesis we have $$Q(N+1,p,\le) \ge N+1- \left \lfloor \frac{N+1}{p+1} \right \rfloor >Q_-(N,p,\le)=Q(N,p,\le).$$ On the other hand, by Property \ref{propQ} (1) we have $Q(N+1,p,\le) \le 1+ Q(N,p,\le)$ so that 
$$Q(N+1,p,\le) = 1+ Q(N,p,\le)=1+N- \left \lfloor \frac{N+1}{p+1} \right \rfloor=Q_-(N+1,p,\le).$$ 
From this last fact and the beginning of the proof we deduce that,  for every $0\le k \le p-1-r_p,$ we have $Q(N_p+k,p,\le)=Q_-(N_p+k,p,\le).$  

By Property \ref{propQ} (4) we may now conclude that for $N \ge N_p$ we have $Q(N,p,\le~)=Q_-(N,p,\le)$ whenever $N \equiv r$ $[p+1]$ for $r_p\le r\le p-1$.

\qed
 
 \bigskip
 Remark that Theorem \ref{prop2} shows only cases where $Q(N,p, \le)=Q_+(N,p,\le).$  In particular Theorem \ref{prop2} proves that $Q(N,p, \le)=Q_+(N,p,\le)$ whenever $p$ is equal to $1$ or $2,$ for any value $N \ge 2p+1$. Hence it is natural to wonder if $Q(N,p, \le)$ is always equal to $Q_+(N,p,\le)$. This is however not the case:
 Wildon~\cite{Wildon16} has checked by a computer search the values of $Q(N,p, \le)$ for $N \le 30$. In particular he listed all such $Q(N,p, \le)$ that are equal to $Q_-(N,p,\le)$. Our (two different) programs  confirmed this list. Looking carefully at the results leads us to two natural questions.

\medskip

\begin{question} \label{N,3,=}
Is it true that for any $N \ge 7$ we have $Q(N,3, \le) = N+1- \lfloor \frac{N+1}{4} \rfloor$~?
\end{question}
\medskip

\begin{question}\label{question2}
Is it true that for any positive integer $p\ge 4$ there exists a smallest integer $N(p) \ge 2p+1$ such that for any integer $N \ge N(p)$ we have: 

$Q(N,p,\le) = Q_-(N,p,\le)$ if and only if $N \equiv r$ $[p+1]$ for $\frac{p+1}{2}< r \le p-1$ ?
\end{question}

\noindent Using the program Main.hs, available from 
Wildon's website\footnote{www.ma.rhul.ac.uk/~uvah099/Programs/MajorityGame/Main.hs},  generating all couples $(N,p)$ such that $Q(N,p,\le) = Q_-(N,p,\le)$ for $N \le 51,$ we could verify that there is no contradiction to a positive answer to Question \ref{question2} and using furthermore Theorem \ref{prop2} and (4) of Property \ref{propQ} we could get that the property is verified for any even $p \le 12$ with $N(4)=9,$ $N(6)=19,$ $N(8)=24$, $N(10)=29$ and  $N(12)=34.$ Notice that, by the proof of Property \ref{propr}, show that for any $p$ there exists always some $N \equiv \lfloor \frac{p+1}{2}\rfloor + 1$ such that $Q(N,p, \le)=Q_-(N,p, \le)$ would be enough to answer positively Question \ref{question2} in case $p$ is even. For the case where $p$ is odd one should furthermore show that $Q(N,p,\le)=Q_+(N,p,\le)$ for every $N \equiv \frac{p+1}{2}$ $[p+1].$

\bigskip

On another hand we could answer (positively) Question \ref{N,3,=}.

\begin{prop} \label{propp?3}
Let $N \ge 7$, we have $$Q(N,3, \le) = N+1- \left \lfloor \frac{N+1}{4} \right \rfloor.$$
\end{prop}

\proof

All along the proof we will use the following property which is similar to Lemma \ref{trouver}.

\begin{property} \label{hauteur}  If a vertex $x$ of a binary binary tree $T_M$ has an associated partition into $u+k$ boxes such that there are more than $2^k$ colorings of the balls that are compatible with the partition, then there exists at least one leaf in the subtree of $T_M$ rooted in $x$ that has at most $u-1$ boxes in its partition.
\end{property}

We know by Theorem 2 that $Q(N,3, \le)$ is equal either to $N- \lfloor \frac{N}{4} \rfloor$ or to 
 $N+1- \lfloor \frac{N+1}{4} \rfloor,$ and that it is equal to $N+1- \lfloor \frac{N+1}{4} \rfloor$ if $N \equiv0,1$ or $3$ $[4].$
 Furthermore, by the computational results we know that Proposition \ref{propp?3} is valid for $N<10$.
 So from now on we may consider that $N=4u+2$ for some $u \ge2$: then $N+1- \lfloor \frac{N+1}{4} \rfloor= N+1-u= N-(u-1)$.

Let $M$ be any method that solves the ($N,3, \le$)-identification problem, we want to prove that the height of the  binary tree $T_M$ is at least $N-(u-1)$, that is, at least one leaf of $T_M$ has an associated  partition into at most $u-1$ boxes.

The tree $T_M$ has a leaf $\ell$ which is connected to the root by a path $P$ whose all edges are labeled with "=". By definition, a partition into boxes is associated to each vertex of $T_M$. We distinguish two cases.

  \begin {itemize}
 \item[Case 1] : No vertex of $P$ has a partition made of boxes each of cardinality either $2$ or $4.$
  
 By definition of $P$ and $\ell$, the partition $\mathcal P(\ell)$ contains only monocolored boxes and it determines  the colors of all the balls, hence all boxes have cardinality at least $4$. Since $N=4u+2$ there are then at most $u$ boxes in  $\mathcal P(\ell)$. In case there are at most $u-1$ boxes in  $\mathcal P(\ell)$ there is nothing more to prove, so let us assume that $\mathcal P(\ell)$ contains exactly $u$ boxes. There are two possibilities for the Type of $\mathcal P(\ell)$ : it is either $(0,4)^{u-2} (0,5)^2$ or $(0,4)^{u-1} (0,6)^1$.
 
 Let $p(\ell)$ be the parent of $\ell$ in $T_M$. If the comparison done in $p(\ell)$ resulted in a ($0,4)$ box, then the partition associated to the brother $v$ of $\ell$ contains $u-1$ monocolored boxes of cardinality at least $4$ and one box of type either $(1,3)$ or $(2,2)$ which leave the colors of some balls undetermined : at least one more comparison will be needed and then there is a vertex whose partition contains at most $u-1$ boxes. The conclusion is similar in case the boxes compared in $p(\ell)$ are of Types $(0,2)$ and $(0,3)$, or both of  Type $(0,3)$. It remains three subcases to consider.
 
 \begin {itemize}
 \item[Subcase 1.1] : The boxes compared in $p(\ell)$ are of Types $(0,1)$ and $(0,4).$ 
 
 Then the Type of $\mathcal P(p(\ell))$ is $(0,1)^1(0,4)^{u-1} (0,5)^1 $ 
 and the Type of
 $\mathcal P(p(p(\ell)))$ is
 
 \begin {itemize}
 \item[(a)] $(0,1)^1(0,2)^2(0,4)^{u-2} (0,5)^1 $  (6 compatible colorings) or 
\item[(b)] $(0,1)^2(0,3)^1(0,4)^{u-2} (0,5)^1 $ (5 compatible colorings) or 
 \item[(c)] $(0,1)^2 (0,4)^{u} $ (4 compatible colorings) or 
 \item[(d)] $(0,1)^1 (0,2)^1(0,3)^1(0,4)^{u-1}$ (5 compatible colorings). \end {itemize}

 Except in case $(c)$ there are at least $5> 2^2$  colorings compatible with the partition into $u+2$ boxes and hence by Property \ref{hauteur} there is a vertex in $T_M$ whose partition contains at most $u-1$ boxes, so the proof is done for Subcase 1.1 (a), (b) and (d). In Subcase 1.1  (c), the parent $p(p(p(\ell)))$ has partition into $u+3$ boxes of Type either $(0,1)^2(0,2)^2(0,4)^{u-1} $ or $(0,1)^3(0,3)^1(0,4)^{u-1} $. In the first case  we have $10$ compatible colorings and in the second case we have $9$, this is more than  $2^3$ and we conclude thanks to Property \ref{hauteur}.

  \item[Subcase 1.2] : The boxes compared in $p(\ell)$ are of Types $(0,1)$ and $(0,5).$
  
  Then the Type of $\mathcal P(p(\ell))$ is $(0,1)^1 (0,4)^{u-1} (0,5)^1 $  exactly as in Subcase 1.1 and we may conclude as well.

 \item[Subcase 1.3] : The boxes compared in $p(\ell)$ are of Types $(0,2)$ and $(0,4).$
 
 Then the Type of $\mathcal P(p(\ell))$ is $(0,2)^1(0,4)^{u}$ 
 and this is not possible by the assumption of Case 1.
 \end {itemize}
 
 \item[Case 2] : There is at least one vertex of $P$ whose partition is made of boxes containing each either $2$ or $4$ balls.
 
 On the path of $T_M$ from the root to $\ell$ there is then a first (nearest from the root) vertex $w$
having the property that in its partition all boxes have cardinality either $4$ or $2$, say the Type of the partition is $(0,2)^{x=2i+1} (0,4)^{y=u-i}$ for some $0 \le i \le u$. Then the boxes that are compared on the parent $p(w)$ of $w$ are both of Type  $(0,1)$ or one is of  Type $(0,1)$ and the other of Type $(0,3)$ and the brother $w'$ of $w$ has a partition of Type either $(0,2)^{x-1} (1,1)^1(0,4)^y $ or $(0,2)^x (1,3)^1(0,4)^{y-1} $. Notice that this partition has exactly one bicolored box, and the small side of this bicolored box has only one ball. Furthermore, each box of $\mathcal P(w')$ contains an even number of balls. Since the partitions of descendants of $w'$ are obtained by merging some boxes, we may conclude that the partitions of the vertices of the subtree of $T_M$ rooted in $w'$ have also all their boxes  containing an even number of balls. 

Then let us consider the leaf $\ell'$ reached from $w'$ using the path $P'$ obtained by allways following an edge labeled "=". On the unique path from the root of $T_M$ to $\ell'$, the vertices on the subpath from the root to $p(w')$ all have a partition into monocolored  boxes and from $w'$ to $\ell'$ the partitions contain only boxes of even cardinality, exactly  one being bicolored, and with only one ball in its small side. Assume the partition associated to $\ell'$ has at least $u$ boxes. By our rule leading to $\ell'$, $\mathcal P(\ell')$ has exactly one bicolored box of Type $(1,z)$. Since there is only one coloring of the balls which is compatible with $\mathcal P(\ell'),$ we have $z \ge 4$ and then all monocolored  boxes should have cardinality at least $3$ and even, hence at least $4$. We may then conclude that the Type of $\mathcal P(\ell')$ is $(1,5)^1(0,4)^{u-1}.$ 

Then,  $\ell'$ cannot be equal to $w'$ and so $\mathcal P(p(\ell'))$ has only boxes of even cardinality and its Type is
 
 \begin {itemize}
 \item[(a)] $(1,5)^1(0,2)^2(0,4)^{u-2}$  and the boxes that are then compared are of Type $(0,2)$, or 
\item[(b)] $(1,3)^1(0,2)^1(0,4)^{u-1}$ and the boxes that are then compared are of Type $(1,3)$ and $(0,2)$ or 
 \item[(c)] $(1,1)^1(0,4)^{u}$ and the boxes that are then compared are of Type $(1,1)$ and $(0,4)$.
 \end {itemize}
 
In the first two cases  the brother of $\ell'$ would have a partition of Type $(1,5)^1(2,2)^1(0,4)^{u-2}$  or $(3,3)^1(0,4)^{u-1}$ that both have at least $2>2^0$ compatible colorings of the balls and we conclude thanks to Property \ref{hauteur}. In case (c), there are two subcases :
\begin {itemize}
 \item[(c1)] $\mathcal P(p(p(\ell')))$ is of Type $(1,1)^1 (0,2)^2 (0,4)^{u-1}$, and then the brother of $p(\ell')$ has partition $(1,1)^1 (2,2)^1(0,4)^{u-1} $ that has  $4>2^1$ compatible colorings and we conclude thanks to Property \ref{hauteur},
 \item[(c2)] $\mathcal P(p(p(\ell')))$ is of Type $(0,1)^2 (0,4)^{u}$ (which means that $p(\ell')=w'$) which is the same as Subcase 1.1 (c).
  \end {itemize}
\end {itemize}
 \qed

\subsection{Bounds and exact values for $Q(N,p, =)$}

In the preceding section we considered the problem of identifying all the colors of the balls for a Red-green coloring with at most $p$ green balls. What happens in  case we know that there are exactly $p$ green balls ?
The next results give upper-bounds of $Q(N,p, =)$. 

\begin{lemma} \label{prop3}
Let $N$ and $p$ be integers such that $1 \le p < \frac{N}{2}$, and let $\nu\ge 0$ be an integer such that $2^{\nu} \le N-2p$, we have :
$$Q(N,p,=) \le 2^\nu -1-\nu +Q(N+1-2^\nu,p,\le).$$

\end{lemma}

\proof  Assume first that $\nu=0$. It is always true that $Q(N,p,=) \le Q(N,p, \le)$, so Lemma \ref{prop3} is verified in the case $\nu=0$.

From now on we will assume that $\nu\ge 1$. As for the preceding proofs, we will use a method which compares boxes, but this time we will have several locations for the boxes : the {\it laboratory} which initially contains all the balls of $\mathcal B$ into boxes of Type $(0,1)$, the {\it reserve} that all along contains only balanced bicolored boxes and a {\it podium} which contains steps numbered starting from $0$: the $i$th step of the podium is provided to receive a monocolored box of size $2^i$.

Any situation where the laboratory contains only boxes of Type $(0,1)$ and the reserve only balanced boxes will be said {\it a correct situation}. The {\it overage of a correct situation} is equal to the overage of red over green among balls that are in the laboratory or in the reserve. Notice that it is the same as the overage among balls in the laboratory since the boxes in the reserve are balanced. 

We now describe a "subroutine" $\mathcal S(i)$, for $i\ge 0$, which starting from any correct situation with an overage $O_i > 2^i$, builds one monocolored box of size $2^i$ , and ends with a correct situation with a new overage of at least $O_i-2^i > 0$. 
\begin{itemize}
\item[-] For $i=0$, any box of the laboratory satisfies the size constraint $2^i$. After putting it aside, the situation is still correct and the overage decreases by at most one (in case the ball in the box we have chosen happened to be red). 
\item[-] For $i=1$, let us compare, one by one, couples of balls in the laboratory: as long as we obtain a balanced box we put it in the reserve.  Since there are at least two more red balls than green balls (the overage is assumed to be at least $2$), we are sure to get once a monocolored couple. Putting it aside, we obtain a correct situation whose overage has decreased by at most $2$.
\item[-]The same way, we may define inductively $\mathcal S(i+1)$ from $\mathcal S(i)$. Assume $\mathcal S(i)$ exists, we will show that then $\mathcal S(i+1)$ exists. 
\newline
Suppose that we are in a correct situation with an overage of at least $2^{i+1}$. We may apply $\mathcal S(i)$ in order to get a monocolored box $B$ of size $2^i$.
We  put $B$ aside. By the induction hypothesis, the situation is then correct and the remaining balls in the laboratory have an overage of at least $2^{i+1}-2^i=2^i$. We may apply $\mathcal S(i)$ again in order to obtain a second monocolored box $B'$ of size $2^i$. We then compare $B$ and $B'$. If this comparison results in a bicolored box, this box is balanced and we put it into the reserve : the overage in the laboratory is still $2^{i+1}$ and we may repeat the procedure. Each time we get  a bicolored box the number of balls in the laboratory decreases, but not the overage, so that we will finally obtain a monocolored box of size $2^{i+1}$. After putting this box aside, by the induction hypothesis, the situation is correct. If the withdrawn balls  are green, the overage of the final situation increases, and else it decreases by at most $2^{i+1}$. So $\mathcal S(i+1)$ does exist.
\end{itemize}

Let $\mu$ be the largest integer such that $2^\mu \le N-2p$ ; since we assumed that $\nu \ge 1$, we have $\mu \ge 1$. 
We will define a method $M(\nu)$ that solves the $(N,p,=)$-identification problem within at most 
$2^\nu -1-m +Q(N+1-2^\nu,p,\le)$ comparisons, thus proving the lemma. We describe now $M(\nu)$.

At Step $0$ of $M(\nu)$ we start with the initial state, that is an empty podium, an empty reserve, and a laboratory containing all the balls in boxes of Type $(0,1)$. So, at this stage, the situation is correct and the overage $\mathcal O(0)=N-2p$ is by definition  at least $2^\mu>\Sigma_{j=0}^{\mu-1}2^j \ge 2^0$.  We apply $\mathcal S(0)$ in order to obtain a monocolored box of cardinality $2^0$ which is put on the $0$th step  of the podium. Then we are in a correct situation with a positive overage $\mathcal O(1)$. If $\mu > 1$ then $\mathcal O(1) > \Sigma_{j=1}^{\mu-1}2^j \ge 2^1$.
We may continue: for each $i \le \nu-1 \le \mu-1$, Step $i$ starts with a correct situation with an overage $\mathcal O(i) >\Sigma_{j=i}^{\mu-1}2^j \ge 2^{i}$ and we apply $\mathcal S(i)$ in order to obtain a monocolored box of cardinality $2^i$ which is put on the $i$th step of the podium. We are then in a correct situation with a positive overage $\mathcal O(i+1)$. If $\mu > i+1 $ then $O(i+1) > \Sigma_{j=i+1}^{\mu-1}2^j \ge 2^{i+1}$.

At the end of Step $\nu-1$ of $M(\nu)$ we have filled $\nu$ steps of the podium with monocolored boxes of cardinality $1, 2, \ldots, 2^{\nu-1}$. Let $\mathcal B_\nu$ be the set of balls that are not on the podium and $N_\nu= \vert \mathcal B_\nu \vert = N-(2^\nu-1)$.  

The reserve contains only balanced boxes and hence an even number, let us say $2k$ ($k \ge 0$), of balls.
The laboratory contains $N_\nu-2k = N+1- 2^\nu-2k$ balls being each contained in one box of Type $(0,1)$ and at most $p-k$ of these balls are green. The overage of red balls in the laboratory is the same as the overage in $\mathcal B_\nu$ and we know that it is positive.

So we may determine the colors of all balls in the laboratory within at most $Q(N_\nu-2k, p-k, \le)$ comparisons and these will lead to a partition $\mathcal P$ of the balls of the laboratory into at least $N_\nu-2k-Q(N_\nu-2k,p-k,\le)$ boxes. By Lemma \ref{borneinfcomm}  we know that $Q(N_\nu-2k,p-k,\le) \le Q(N_\nu,p,\le)-2k$ ; from this we deduce that $\mathcal P$ has at least $N_\nu-Q(N_\nu,p,\le)$ boxes. Comparing one of these boxes with all boxes in the reserve will then provide a partition $\mathcal P_\nu$ of $\mathcal B_\nu$ which determines completely the colors of the balls in $\mathcal B_\nu$ and has the same number of boxes than $\mathcal P$. 

Then we know exactly the number $p' \le p$ of green balls in $\mathcal B_\nu$ and there should be $p-p'$ green balls on the podium. The binary representation of $p-p'$  gives us the boxes of the podium containing the green balls without any further comparisons. So we have a  partition of the whole set of balls  to which corresponds a unique Red-green coloring and it consists into $\nu$ boxes for balls on the podium and at least $N_\nu-Q(N_\nu,p,\le)$ boxes for balls in $\mathcal B_\nu$. 
The number of comparisons that are done to obtain such a partition is then at most  $$N-(\nu+N_\nu-Q(N_\nu,p,\le)=2^\nu-1-\nu+Q(N-(2^\nu-1),p,\le).$$
\qed

\medskip

As a consequence of the preceding results we obtain the following theorem.

\begin{theorem}\label{th=}
Let $N$ and $p$ be integers such that $1 \le p < \frac{N}{2}$, and let $m$ be the largest integer such that $2^{m} \le \min(N-2p, 2p)$, we have :
$$Q(N,p,=) \le N+1-m - \left \lfloor \frac{N+2-2^m}{p+1} \right \rfloor.$$
\end{theorem}

\proof
Combining Lemma \ref{prop3} and Theorem \ref{prop2} we get that $$Q(N,p,=) \le \min_{0 \le \nu \le \mu} \left (N+1-\nu - \left \lfloor \frac{N+2-2^\nu}{p+1} \right \rfloor \right).$$ It remains to show which value of $\nu$ provides the minimum of $f(x)= N+1-x - \lfloor \frac{N+2-2^x}{p+1}\rfloor.$

We will do it by comparing the value of $f$ for two consecutive values $\nu$ and $\nu+1$ such that $2^{\nu+1} \le N-2p$. We may write $N= 2p+2^{\nu+1}+ \ell=2p+2^\nu+2^\nu + \ell$ where $\ell \ge 0$:

\begin{itemize}
\item $f(\nu) = N+1-\nu -\lfloor \frac{N+2-2^\nu}{p+1}\rfloor = N+1-\nu -\lfloor \frac{2p+2^\nu + \ell+2}{p+1}\rfloor=N-\nu-1-\lfloor \frac{2^\nu + \ell}{p+1}\rfloor$
\item $f(\nu+1) = N-\nu -\lfloor \frac{N+2-2^{\nu+1}}{p+1}\rfloor = N-\nu -\lfloor \frac{2p+\ell+2}{p+1}\rfloor = N-\nu -2-\lfloor \frac{\ell}{p+1}\rfloor $
\end{itemize}

So, if $2^\nu\le p$ then $f(\nu+1) \le f(\nu)$ and else $f(\nu+1)\ge f(\nu)$. 

In case $N-2p \le 2p$, then any $\nu+1\le \mu$ verifies that $2^{\nu+1} \le 2^\mu \le N-2p \le 2p$ and then $2^\nu \le p$. So, then
$\min_{0 \le \nu \le \mu}(N+1-\nu - \lfloor \frac{N+2-2^\nu}{p+1}\rfloor)$ is attained for $\nu= \mu$, so $m=\mu$. 

If $2p< N-2p$, then the minimum will be attained for $m$ equal to the largest $\nu$ such that $2^{\nu-1} \le p$, or equivalently $2^\nu \le 2p$.
\qed

\medskip

\begin{prop} \label{propp=1}
Let $N \ge 3$,  we have :
$$Q(N,1,=) = N- \left \lfloor \frac{N}{2} \right \rfloor.$$
 \end{prop}

\proof 
For $N=3$ it is easy to verify that the proposition holds.
For $N > 3$, Theorem \ref{th=}  gives directly $Q(N,1,=) \le N- \lfloor \frac{N}{2}\rfloor$.
Let us show now that we also have $Q(N,1,=) \ge N- \lfloor \frac{N}{2}\rfloor$.

Let $M$ be any method that solves the $(N,1, =)$-identification problem, $T_M$ be its associated binary tree, and $P=(R=u_0,\ldots, u_k=l)$ be the path in $T_M$, from the root to a leaf $l$, whose all edges are labeled "$=$". The partition $\mathcal P(l)$ consists in monocolored boxes, and exactly  one $1$-equal Red-green coloring should be compatible with $\mathcal P(l)$. Hence it contains exactly one box of cardinality $1$, all others being of cardinality at least $2$. 

If there are at most $\lfloor \frac{N}{2}\rfloor$ boxes in $\mathcal P(l)$ then the height of $T_M$ is at least $N- \lfloor \frac{N}{2}\rfloor$ and  the desired inequality is verified. 

When $N$ is even, a partition of the balls verifying the property above cannot have more than $\frac{N}{2}=\lfloor \frac{N}{2}\rfloor$ boxes (this bound is reached if and only if one box of $\mathcal P(l)$ has cardinality $1$, another one has cardinality $3$ and all others have cardinality $2$). It remains to consider the case where $N$ is odd and $\mathcal P(l)$ has at least $\lceil \frac{N}{2}\rceil$ boxes. This can happen only if all boxes of $\mathcal P(l)$ except one have cardinality $2$. Consider now the brother $l'$ of $l$, that is the other child of the parent $u$ of $l$. The partition of $l'$ is the same as the one of $l$ except that there is exactly one bicolored box $B$, of Type $(1,1)$. Then the green ball can be either of the two balls in $B$, and one more comparison should be done in order to conclude. Hence the height of $T_M$ is at least $N- \lfloor \frac{N}{2}\rfloor$, and the proof is done.
\qed

\begin{prop} \label{propp=2}
Let $N \ge 6$,  we have :
$$Q(N,2,=) = N- \left \lfloor \frac{N+1}{3} \right \rfloor.$$
 \end{prop}

\proof 
We have verified by a computer program that this proposition holds for $N \le 9$, so we may assume that $N\ge 10$.
Furthermore, for any $N \ge 8$, the value of $m$ as defined in Lemma \ref{prop3} is $2$, and so $Q(N, 2, =) \le N + 1 - 2 - \lfloor \frac {N+2-2^2}
{2+1} \rfloor = N - \lfloor \frac {N+1} {3} \rfloor.$
Then  it is sufficient to show that for any $N \ge10$: $$Q(N, 2, =) \ge  N - \left \lfloor \frac {N+1} {3} \right \rfloor.$$ 

Let $M$ be any method that solves the $(N,2, =)$-identification problem, $T_M$ be its associated binary tree, and $P=(R=u_0,\ldots, u_k=l)$ be the path in $T_M$, from the root to a leaf $l$, whose all edges are labeled "$=$". For each vertex of $P$, the associated partition consists only of monocolored boxes. A ball  will be said {\it meager} if it is in a monocolored box of cardinality $1$ or $2$, else it is called {\it fat}. In the partition associated to the root, all $N$ balls are meager and in $\mathcal P(l)$ at most $3$ balls are meager (it is easy to verify that else there would be at least two $2$-equal Red-green colorings of the balls compatible with $\mathcal P(l)$). After each comparison the number of meager balls changes as follows :

\begin{itemize}
\item[-] it remains the same after comparing two fat balls or two balls in cardinality~$1$ boxes,
\item[-] it decreases by $1$ after comparing a fat ball and a ball in a cardinality $1$ box,
\item[-] it decreases by $2$ after comparing a fat ball and a ball in a cardinality $2$ box,
\item[-] it decreases by $3$ after comparing a ball in a cardinality $1$ box and a ball in a cardinality $2$ box,
\item[-] it decreases by $4$ after comparing balls in cardinality $2$ boxes.
\end{itemize}

From the remarks above we may define $k$ as the largest index such that  $\mathcal P(u_k)$ has at least $10$ meager balls. Then by definition of $k$ and the possible decreases after one comparison, $\mathcal P(u_{k+1})$ contains at most $9$ and at least $6$ meager balls.

Let us assume that $\mathcal P(u_{k+1})$ contains $n_1$ boxes of Type $(0,1)$, $n_2$ boxes of Type $(0,2)$ and $n_3$ monocolored boxes of cardinality at least $3$. We have $k+1= N-(n_1+n_2+n_3)$,  $m=n_1+2n_2$ is the number of meager balls of  $\mathcal P(u_{k+1})$ and $6 \le m \le 9$.

Obviously, in $\mathcal P(u_{k+1})$,  the two green balls are either in one box of cardinality $2$ or in two boxes of cardinality $1$. This means that there are exactly ${n_1 \choose 2} + {n_2 \choose 1}$ $2$-equal Red-green colorings compatible with $\mathcal P(u_{k+1})$ (we set ${n_1 \choose 2}=0$ for $n_1 \le 1$). 
Then the subtree of $T_M$ induced by $u_{k+1}$ and all its descendants has ${n_1 \choose 2} + {n_2 \choose 1}$ leaves. By Lemma \ref{trouver} we know then that 
  at least $q(n_1,n_2)=\lceil \log_2 ({n_1 \choose 2} + {n_2 \choose 1}) \rceil$ additional comparisons are necessary in order to have  this number of leaves in the subtree. So, the height of $T_M$ is at least  $k+1+q(n_1,n_2)$.  To end the proof, it is enough to show that 
$$k+1+q(n_1,n_2) \ge N- \left \lfloor \frac{N+1}{3} \right \rfloor$$
or equivalently (since $k+1+q(n_1,n_2)$ is an integer) $$3N-3n_1-3n_2-3n_3+3q(n_1,n_2) \ge 3N - N-1.$$
Since $N = n_1+2n_2+3n_3+e$ for some $e\ge 0$ this is the same as $$e+1 \ge 2n_1+n_2-3q(n_1,n_2).$$ As $6 \le m \le 9$ and $n_1$ and $n_2$ should be such that $n_1+2n_2=m$, there is only a finite number of cases to consider. These are summarized in Table \ref{e+1}. Since in all cases we get $2n_1+n_2-3q(n_1,n_2) \le 1$, the needed inequality indeed holds.
\qed

\begin{center}
\begin{table}[ht!] 
{\footnotesize \begin{tabular}{|m{1,5cm}||c|c|c|c|c|c|c|c|c|c|c|c|c|c|c|c|c|c|}
\hline
$m$ & 6 & 6 & 6 & 6 & 7 & 7 & 7 & 7 & 8 & 8 & 8 & 8 & 8 & 9 & 9 & 9 & 9 & 9 \\
\hline
$n_1$ & 6 & 4 & 2 & 0 & 7 & 5 & 3 & 1 & 8 & 6 & 4 & 2 & 0 & 9 & 7 & 5 & 3 & 1\\
\hline
$n_2$ & 0 & 1 & 2 & 3 & 0 & 1 & 2 & 3 & 0 & 1 & 2 & 3 & 4 & 0 & 1 & 2 & 3 & 4\\
\hline
\hline
$2n_1+n_2-3q(n_1,n_2)$ & 0 & 0 & 0 & -3 & -1 & -1 & -1 & -1 & 1 & 1 & 1 & 1 & -2 & 0 & 0 & 0 & 0 & 0\\
\hline
\end{tabular}}
\caption{The value of $2n_1+n_2-3q(n_1,n_2)$ for $6 \le n_1+2n_2=m \le 9$.}\label{e+1}
\end{table}
\end{center}

Propositions \ref{propp=1} and \ref{propp=2} show that, for $p=1$ or $2$, the upper bound $Q_+(N,p, =)$ of $Q(N,p, =)$ given by Theorem \ref{th=},
is in fact the right value of $Q(N,p,=)$. This is however not true for any value of $p$ as shown by our programs.
For $N \le 30$ and $p < \frac{N}{2}$ it happens several times that $Q(N,p,=)=Q_+(N,p,=)-1$ (see Table \ref{diff1}),  once that $Q(N,p,=)=Q_+(N,p,=)-2$ (for $(N,p)= (30,4)$)
and for all other couples $(N,p)$ the equality $Q(N,p,=)=Q_+(N,p,=)$ holds.

\begin{center}
\begin{table}[ht!] 
{\footnotesize \begin{tabular}{|p{0,05cm}||c|c|c|c|c|c|c|c|c|m{0,7cm}|m{1,0cm}|c|}
\hline
$N$ & 17 & 20 & 21 & 22 & 23 & 24 & 25 & 26 & 27 & 28 & 29 & 30 \\
\hline
$p$ & 3 & 3, 4 & 3 & 5 & 5 & 3, 4 & 3, 4, 6 & 6 & 3, 4, 5 & 3, 4, 5, 7 & 3, 4, 5, 6, 7  & 3, 6, 8 \\
\hline
\end{tabular}}
\caption{The values of $N$ and $p$ such that $Q(N,p,=)=Q_+(N,p,=)-1$ .}\label{diff1}
\end{table}
\end{center}

From these results one could expect that  $Q(N,p,=)$ is never very far from its upper bound $Q_+(N,p,=)$. By the results below, this is however not the case as soon as $p$ is at least $3$. 

In this part it will be convenient to use the notation $B(N,p,=)$ for the value $N-Q(N,p,=)$. This value represents the minimum number of boxes we may have at the end of any optimal method solving the $(N,p,=)$-problem (i.e. using at most $Q(N,p,=)$ comparisons). Similarly we define $B_+(N,p,=)=N-Q_+(N,p,=).$ (Notice that $B_+(N,p,=)$ is then a lower bound of $B(N,p,=)$). By Theorem \ref{th=} we know that for a "sufficiently big $N$", $B_+(N,p,=)$ is approximately $\frac{N}{p+1}$. We will show in the following that a better ratio may be obtained for any $p\ge 3$. Let us first consider the case where $p=3$.

\begin{prop} \label{propp=3}
Let $N \ge 7$,  we have for some positive constant $D$:
$$B(N,3,=) \ge  \left \lfloor \frac{3N}{10}\right \rfloor - D.$$
 \end{prop}

\proof
To prove the bound we exhibit a method giving the colors of all balls that ends with at least $\lfloor \frac{3N}{10}\rfloor - D$ boxes, for a constant $D>0$.

We first compare $b_1$ to $b_2, b_3, ...$ until to obtain $4$ balls of the same color (and hence red). Thus we get a partition containing exactly one bicolored box and, more precisely, of Type $(u,4)^1 (0,1)^{N-(u+4)}$ for some integer $u \le 3.$ 

Then we build monocolored boxes of cardinality alternately $2$ and $3$ : during this process, each time we get a bicolored box we compare it to the previous lonely bicolored box, this will happen at most $3-u$ times, and doing so we keep the property of having exactly one bicolored box.
We stop this process when there is no more cardinality 1 boxes. 
At this stage, if we don't have the same number of boxes of cardinality $2$ and $3$, we compare the boxes obtained during the last trial with the unique bicolored box, increasing its cardinality by at most $4$. 
So at the end of this process, we have a partition of the balls of Type $(x,y)^1 (0,2)^a(0,3)^a$
where $y\ge 4$, $x+y \le 17$ and $ a\ge \frac{N-17}{5}$. It is not possible that $x=2$, since else the remaining green ball cannot be in a box of Type $(0,2)$ or $(0,3)$. So either $x=3$ and we know the colors of all balls, or $x=0$ and the three green balls are in a monocolored box of size $3$ or $x=1$ and the two remaining green balls are in a monocolored box of size $2$. In the last two cases, we have to solve a problem equivalent to the one of finding one green ball among $a$ balls, which by Proposition \ref{propp=1} can be done leaving at least $\frac{a}{2}$ boxes. So in total we have at least $1+a+\frac{a}{2}$ boxes which is  at least $1+\frac{3a}{2}\ge\frac {3N}{10}-4.1$, so more than $\frac{3N}{10} -5$. \qed

We observe that the ratio $\frac{3}{10}$ is slightly better than the ratio  $\frac{1}{4}$ obtained by Theorem~\ref{th=}. Hence, for a big enough $N$ the method described in the proof of Proposition \ref{propp=3} is more efficient than the one in the proof of Theorem~\ref{th=}.

\subsection{The Towers Method}

The method used in the proof of Proposition \ref{propp=3} can be generalized for any value of $p$ by the {\it  Towers Method} that we describe now : this name comes from the fact that we may consider that we build, from the boxes obtained by the chosen comparisons, towers each containing all monocolored boxes of a given size, the {\it height of a tower} corresponding to the number of boxes it is made of. For a given $p$ the method depends on the following parameters :
\begin{itemize}

\item a non empty set $E=\{u_1, \ldots,u_k\} \subseteq \{2, \ldots,p\}$ of {\it sizes of the monocolored boxes} that will be built. This subset should satisfy the following "{\bf unicity property}" :
For every $0\le l \le p,$ there exists at most one $k$-tuple of non-negative integers $(q_1, \ldots, q_k)$  such that $l=q_1u_1 + \ldots +q_ku_k$. We denote by $\mathcal U(p)$ the set of subsets of $\{2, \ldots,p\}$ that have the unicity property and by $L_E$, $E$ being in $\mathcal U(p)$, the set of values $0\le l \le p,$ for which there exists a $k$-tuple associated with the above expression of $l$ as a sum of values in $E.$

The cardinality $k \le p-1$ of $E$ will be the number of towers that will be erected. In the proof of Proposition \ref{propp=3} we had
$E=\{u_1,u_2\}$ with $u_1=2$ and $u_2=3$ and indeed $1$ cannot be expressed as a sum of $2$'s and $3$'s, and $0$, $2$ or $3$ can, but by only one way ; so $L_E=\{0,2,3\}.$

\item a set of $k$ integers $h_1,  \ldots h_k,$ expressing each the number of monocolored boxes of size $u_i$ that are added to the $u_i$-tower at each step. In the proof of Proposition \ref{propp=3} we had $h_1=h_2=1$.

\end{itemize}

Now we may describe the Towers method. 

As usual all balls are initially in boxes of Type $(0,1)$. At the first step we compare $b_1$ to $b_2, b_3, ...$ until to obtain $p+1$ balls of the same color (and hence red). Thus we get a partition containing exactly one bicolored box of Type $(u,p+1)$ for some integer $u \le p,$ with all other boxes of cardinality $1.$

At the second step we use a procedure "add floors to each tower" by  consecutively building $h_1$ monocolored boxes of cardinality $u_1$ added to the $u_1$-tower, $h_2$ monocolored boxes of cardinality $u_2$ added to the $u_2$-tower, \ldots, $h_k$ monocolored boxes of cardinality $u_k$ added to the $u_k$-tower. During this process, each time it happens that we create a new bicolored box we compare it to the previous lonely bicolored box in order to keep a partition with exactly one bicolored box. We continue to "add floors" as long as we can. If we have to stop because of lack of cardinality 1 boxes before finishing the "add floors" procedure, then we compare  each box created during the last and partial "add floors" procedure with the unique bicolored box. At that point, the partition is of Type $(x,y)^1(0,u_1)^{h_1c} \ldots(0,u_k)^{h_kc}$ where $x\le p$, $y\ge p+1$, $c=\frac{N-(x+y)}{h_1u_1+\ldots +h_ku_k}$ (so $c$ corresponds to the number of times we could perform entirely the procedure "add floors"). 

At this stage we know for sure $x$ green balls and we have to detect the remaining $l=p-x$  ($0 \le l \le p$). Since these balls are in the towers, by the "unicity property", there exists a unique $(q_1(l), \ldots, q_k(l))$  such that $l= \Sigma_{1 \le i \le k}q_i(l)u_i $ and $q_i(l) \ge 0$ for each $i$ ; that is $l \in L_E$.
It remains now to solve the $(h_ic,q_i(l),=)$-problem for each $1\le i \le k$ in order to find all green balls, since one monocolored box may be considered as just one ball of the same color. Any method solving these problems may be used. In case we know for each $(N,q_i(l),=)$-problem a constant $C(q_i(l))$ such that the $(N,q_i(l),=)$-problem may be solved leaving at least $C(q_i(l))N -D_i$ boxes, for some constant $D_i$, since
exactly as in the case of the proof of Proposition \ref{propp=3} the value of $x+y$ is bounded by a constant, we would then get at least  $\frac {\sum_{1 \le i\le k} C(q_i(l))h_i}{h_1u_1+\ldots +h_ku_k}N-D_l$ remaining boxes, for some constant $D_l$. Notice that we know from Theorem \ref{th=} that coefficients $C(q_i(l))\ge \frac{1}{q_i(l)+1}$ do exist.

As we cannot fix the value $l \in L_E$ of balls remaining to be detected we obtain
$B(N,p,=) \ge \min_{l\in L_E} \{\frac {\sum_{1 \le i\le k} C(q_i(l))h_i}{h_1u_1+\ldots +h_ku_k}\}N-D,$ for some constant $D$. 

Notice that since $1 \notin E$ we have $q_i(l) <p$. Hence 
we may show now how, based on best known values $C^T(i)$ for $i <p$ we may compute recursively, for each $p \ge 3,$ the Towers Method's parameters providing the maximum constant  $C^T(p)$, by solving a linear programming problem. We will need values of $C^T(i)$ for $i<3$: the $(N,0,=)$-problem can obviously be solved leaving $N$ boxes. By Propositions \ref{propp=1} and \ref{propp=2}, $B(N,1,=)\ge \frac{1}{2}N-\frac{1}{2}$ and $B(N,2,=)\ge \frac{1}{3}N-\frac{1}{3}$. So we set $C^T(0)=1$, $C^T(1)=\frac{1}{2}$ and $C^T(2)=\frac{1}{3}$.

First we have to compute all subsets $E=\{u_1, \ldots,u_k\}\in \mathcal U(p)$. 
For each such $E$, we will show how to compute heights $h_1, h_2, \ldots h_k$ providing the best proportion $C^E(p)$ of boxes we are sure to obtain by a Towers Method based on $E.$ For that purpose, we need to determine what is in such an optimal method, the proportion $x_r$ of balls placed in monocolored boxes of size  $u_r$ among the set of  balls placed in monocolored boxes, for each $1\le r \le k$. So $x_r $ represents $\frac{h_ru_r}{h_1u_1+\ldots +h_ku_k}$. We will need one more variable $y_E$, and the linear programming problem $P_E$ will be the following :

Maximize $y_E$ subject to
$\begin{cases} 
\sum_{i=1}^{k}x_i=1, \\
y_E - \sum_{1 \le i\le k} C(q_i(l))u_i^{-1}x_i \le 0  &  \text{for each } l \in L_E,  \\
x_i \ge 0 \\
\end{cases}$

It is clear that $P_E$ is feasible. Let $(x_1,x_2, \ldots,x_k,y_E)$ be an optimal vertex of $P_E.$ Since all coefficients of $P_E$ are rational we have that all $x_i$'s are rational numbers, and then there exist integers $a_1, \ldots, a_k,b$ such that $x_r=\frac{a_r}{b}= \frac{a_r \Pi_{1\le i \le k}u_i}{b\Pi_{1\le i \le k}u_i}$ for each $1 \le r \le k$. A Towers Method with parameters $E$ and $h_r=a_r \Pi_{1\le i \le k, i\neq r}u_i$ for each $1\le r\le k,$ will provide, for a big enough number of balls, a proportion $x_r$ of balls placed in monocolored boxes of size $u_r$ among balls placed in monocolored boxes.

We can solve $P_E$ for each $E \in \mathcal U(p)$ and the best ratio $C^T(p)=\max_{E \in \mathcal U(p)}y_E.$ Notice that it is enough to consider only maximal subsets $E$ in $\mathcal U(p)$.

In particular one can show with this method that $C^T(3)=\frac{3}{10}$, $C^T(4)=\frac{5}{18}$, $C^T(5)=\frac{2}{9}$ and $C^T(6)=\frac{1}{5}$.

However one may also wonder if there are other methods that could give better bounds.

For any positive integer $p$, let us define $I_p$ as the set of values $C\in [0,1]$ for which there exists a constant $D \ge 0$ such that, for each $N> 2p,$ we have  $$B(N,p,=) \ge CN - D.$$
Let $C_p$ be the supremum of $I_p$. It is easy to see that $I_p$ is equal either to $[0,C_p]$ or to $[0,C_p[,$ and that, by Theorem \ref{th=}, $C_p \ge \frac{1}{p+1}$ for every $p\ge0$.

\begin{prop} \label{limc}
The sequence $(C_p)$ is non-increasing and $\lim\limits_{p \rightarrow \infty} pC_p = \infty.$
\end{prop}

\proof Let us consider an integer $p\ge 1$ and an integer $N \ge 2p+1$. By Lemma \ref{borneinfcomm} one has $$B(N-2,p-1,=)=N-2-Q(N-2,p-1,=) \ge N-Q(N,p,=)= B(N,p,=).$$ Assume that for some $C>0$ there exists $D \ge 0$ such that  $B(N,p,=)\ge CN-D.$ Then by the previous inequality we get that $B(N-2,p-1,=)\ge C(N-2)-D$ and hence $C_{p-1} \ge C_p.$

To state $\lim\limits_{p \rightarrow \infty} pC_p= \infty,$ we will use a coefficient provided by the Towers Method.

Given an integer $q\ge 2$,  let $p$ be any integer such that $p \ge q^q$ and set $k_p=\lfloor \frac{p}{q}\rfloor$.

\begin{claim} \label{unicity}
The set $E=\{u_i \vert u_i = k_p + q^{i-1}, i\in\{1,\ldots,q-1\}\}$ belongs to $\mathcal U(p).$
\end{claim} 

In order to prove this Claim we have to show that for any $0 \le l \le p$ there exists at most one $(q-1)$-tuple $(x_1, \ldots,x_{q-1})$ of non-negative integers such that $l=\Sigma_{1 \le i \le q-1} x_iu_i.$
So, let us consider $0 \le l \le p$ and let $s=\lfloor\frac{l}{k_p+1}\rfloor.$ Since $l \le p,$ by the definition of $k_p,$ we have $0\le s \le q-1$. Then, since each $u_i \ge k_p+1$ we have $\Sigma_{1 \le i \le q-1} x_i \le s.$

On another hand we have $k_p \ge q^{q-1} > (q-2)q^{q-2} \ge (s-1)q^{q-2}-s$, and then by definition of $s$ we have $l \ge s(k_p+1)=(s-1)k_p+k_p+s > (s-1) (k_p + q^{q-2}).$ Then we have $\Sigma_{1 \le i \le q-1} x_i \ge s.$ By the inequality above we conclude that $\Sigma_{1 \le i \le q-1} x_i = s.$ There is a unique way to write $l-sk_p$ in basis $q$ and hence at most one way to obtain $l$ balls as the union of boxes of cardinalities $k_p+q^0, \ldots, k_p+q^{q-2}.$ The Claim is proved.
 
Let us now assume that for some "big" $N$ we use the Towers method with $E$ as defined in Claim \ref{unicity} and heights $h_i=1$ for each $1 \le i \le q-1.$ Let $l$ be the number of green balls that should be detected at the end of the erection of the towers. Notice that the towers all have the same number $a$ of monocolored boxes. By the proof of Claim \ref{unicity} there are $s=\lfloor\frac{l}{k_p+1}\rfloor$ monocolored boxes to discover and $0\le s \le q-1$. We claim that the case leaving the smallest number of boxes is when there is one green box in each of the $q-1$ towers, in which case we end with at least $(q-1)\frac{a}{2}$ boxes. Indeed if we have less than $q-1$ green boxes to find, then obviously we won't end with less boxes, and in case we have to find $q-1$ green boxes among at most $q-2$ towers we gain at least $\frac{a}{2}$ boxes and loose at most $\frac{a}{2}-C_{q-1}a \le \frac{a}{2}$. 
It remains to observe that $N$ is approximately equal to $$a\Sigma_{1\le l \le q-1} (k_p + q^{l-1}) = a \left((q-1)k_p+ \frac{q^{q-1}-1}{q-1}\right) \le aqk_p \le ap.$$

Hence $B(N,p,=)\ge (q-1) \frac{N}{2p}-R$ for some constant $R$ which means that $C_p \ge \frac{q-1}{2p}$ and then $pC_p \ge \frac{q-1}{2}$, thus $\lim\limits_{p \rightarrow \infty} pC_p = \infty.$ \qed

\section{Some additional remarks}

\subsection {Other formulations}
The problems we considered here on a set of bicolored balls are also known with other equivalent formulations. 

One example is : we are given a set of coins
having two possible weights, some coins are faked, we don't know if they are heavier or lighter, but we know that there are less fake coins than true coins and we may only compare the weights of two coins with a two-pan balance. 

Another formulation is the following. The population of an island is composed of "knights", who always tell the truth, and "knaves", who always lie (these names come from Smullyan book \cite{Smu} but the same problem appears with several other terms). There are more knights than knaves in the population and everyone knows the status of each. A visitor is allowed to choose couples $(a,b)$ of inhabitants of the island and ask $a$ : "Is $b$ a knight ?". A"yes" answer to the question will always mean that $a$ and $b$ are of the same kind and a "no" answer  that they are of different kinds.  So this is exactly the result of a comparison... 

In any of these formulations one may be interested by finding one object in the majority, or one object in the minority, or the status of all objects, ..., by doing a minimum number of comparisons.
To our knowledge the case where the cardinalities of the subsets are fixed has not been studied before.

The problem above is a particular case of a more general problem, where it is assumed that the inhabitants are divided into "knights" and "spies", with a majority of knights. A spy can either lie or tell the truth.
This problem too has been studied with other names. In particular there is a formulation in terms of chips testing \cite{AlonsoEtAl04}. Notice that surprisingly the minimum number of questions needed to find a knight is the same in the "knights and knaves" and in the "knights and spies" cases \cite{AlonsoEtAl04}: $N - b(N)$ (as explained at the beginning of the paper for the Majority problem). In fact Alonso et al \cite{AlonsoEtAl04} showed that any algorithm which solves the Majority problem can be used to solve this knights and knaves problem, and both algorithms behave the same in the worst case. However this seems not to be true for the average-case \cite{AloReiSch97}. Also, the minimum number of questions to ask in order to know the status of everyone is different: $N+p-1$ \cite{Ble}.

\subsection {A game}
The $(N,p,\le)$- and $(N,p,=)$-identification problems can be viewed as a game between two players named "Maker"  and "Breaker". At the beginning of each step,  Maker chooses the two balls to compare and then  Breaker decides the answer of the comparison. The answers however, should always allow at least one compatible $p$-majored (resp. $p$-equal) Red-green coloring of the balls. 

The two players agree on an integer $0 \le q \le N-1$. The winner is  Maker if, after at most $q$ steps, he is able to know the colors of all balls, and else it is Breaker. So, clearly,  Maker is sure to win as soon as $q \ge Q(N,p, \le)$ (resp. $q \ge Q(N,p, =)$). 

The proofs of the lower bound for the Majority problem by Saks and Werman \cite{SakWer91}, and later \cite{Wie}, use this game approach. 
 Blecher \cite{Ble} used it also for the Knights and Knaves problem described above. Other authors \cite{Aig,Wildon10} used this game approach for similar problems.

\subsection {The ultimate formulation}
The most abstract formulation of a game corresponding to the $(N,p,=)$-identification problem is the following.

The game begin with two integers $N$, $p$  such that $0 \le 2p < N$, and an initial pattern of $N$ couples $(1,0)$.

After $N-k$ steps we will have a pattern of  $k$ unordered couples 
$(u_1,v_1) \ldots (u_k,v_k)$  such that
the $u_i$'s and  $v_i$'s are non negative integers. 

We define $C((u_1,v_1) \ldots (u_k,v_k), p,=)$ as being equal to the value of :

\noindent $\vert \{ (r_1, \ldots, r_k) \in \Pi_{1 \le i \le k} \{u_i,v_i\} \vert r_1+ \ldots + r_k=p\}\vert.$ 

The next pattern is obtained as follows : Maker chooses two indices $l$ and $m$ such that $1\le l < m \le k$,
 we will keep in the pattern all $(u_i,v_i)$'s, $l \neq i \neq m$,
and we will merge the couples $(u_l,v_l)$ and $(u_m,v_m)$ replacing them by a unique couple
which will be either $(u_l + v_m, v_l + u_m)$ or $(u_l + u_m, v_l + v_m)$. It is  
 Breaker who decide which one of these couples
has to be chosen, with the additional condition that, denoting the new pattern by $(u'_1, v'_1)\ldots(u'_{k-1}; v'_{k-1})$,
we have $C((u'_1, v'_1)\ldots(u'_{k-1}, v'_{k-1})) > 0$.

Maker and Breaker agree on an integer $q$ (with $0 \le q\le N-1$) and Breaker wins if $C((u_1, v_1), \ldots, (u_{N-q}, v_{N-q)}, p, =) > 1$,
otherwise Maker wins.

Clearly the minimum value of $q$ allowing Maker (assuming that he plays perfectly)
to be sure to win, whatever 
Breaker plays, is $Q(n, p, =)$.

Indeed this is the same problem as the bicolored balls,
but we limit ourselves to a minimal quantity of information (the one given by the patterns)
without even speaking of the set of $N$ bicolored balls.

We remark that $C((u_1, v_1), \ldots,  (u_k, v_k), p, =)$ is nothing else than the number of possible colorings  with the pattern
$(u_1,v_1) \ldots (u_k, v_k)$.

Without loss of generality we may assume that we merge $(u_1, v_1)$ and $(u_2 v_2)$ and we check easily that:
$C((u_1+u_2, v_1+v_2)(u_3, v_3),\dots, (u_k, v_k), p, =)+C((u_1+v_2, v_1+u_2)(u_3, v_3), \ldots, (u_k, v_k), p, =)$
is equal to $C((u_1,v_1), \ldots, (u_k,v_k), p,=)$.

Breaker may use a "majority strategy":
when Maker has chosen to merge the couples $(u_l, v_l)$ and $(u_m, v_m)$, Breaker will choose the couple among 
$(u_l + v_m, v_l + u_m)$ and $(u_l + u_m, v_l + v_m)$ which maximize 
the number of  colorings of the subsequent patterns (or any of the couples if this number of coloration is the same).

\begin{rema} If we have the pattern
$(u_1,v_1)\ldots(u_k,v_k)$ and if  Breaker uses a "majority strategy"
then the preceding relation shows that Maker will need at least 
$\lceil \log_2 C((u_1, v_1), \ldots, (u_k,v_k))\rceil$ more steps to get a number of colorings equal to $1$.
\end{rema}

We also notice that this strategy allows Breaker to be sure that 

\noindent $C((u'_1, v'_1) \ldots (u'_{k-1}, v'_{k-1})) > 0$ in case we had already 

\noindent $C((u_1, v_1) \ldots (u_k, v_k)) > 0$.

On his side, Maker may apply a strategy (for the choice of $l$ and $m$) which minimizes the number of colorings in case Breaker uses a "majority strategy".

What has just been done concerns the case of the $(N,p, =)$-identification problem
but it is easy to find an analogous interpretation for the case of the $(N,p, \le)$-identification problem.
We remark that the found upper bounds for $Q(N,p,=)$ and $Q(N,p,\le)$
are the work of Maker and the found lower bounds for $Q(N,p,\le)$ and $Q(N,2,=)$ are the work of Breaker.

\bigskip
Question : Is it possible to improve the strategies of respectively Maker and Breaker in
 a way to give better  bounds for $Q(N, p, \le)$
and $Q(N, p, =)$ ?
\bigskip

Complementary remarks :

1) In the patterns, we may use the multiplicative notation, replacing 
$r$ times the couple $(u, v)$ : $(u, v)\ldots(u, v)$, 
by $(u,v)^r$.

2) Maker may decide to merge more than two couples regardless of the answers of Breaker.
Example : if $N = 17$, Maker may (at the beginning of the game) merge
four times three couples $(1, 0)$, the result will be a pattern $(1, 0)^5$ $(3, 0)^k$ $(2, 1)^{4-k}$, for $0\le k \le 4$.

3) The initial problem could be generalized for more than two colors.

\bigskip

\bibliographystyle{plain}
\bibliography{biblio}

\end{document}